# Hexagonal Surfaces of Kapouleas

Frank Morgan
Department of Mathematics and Statistics
Williams College
Williamstown, Massachusetts 01267
Frank.Morgan@williams.edu

## Abstract

For the "hexagonal" norm on $\mathbf{R}^3$, for which the isoperimetric shape is a hexagonal prism rather than a round ball, we give analogs of the compact immersed constant-mean-curvature surfaces of Kapouleas.

## 1. Introduction

In 1986, H. Wente [W] provided the first example, other than a round sphere, of a *compact* (immersed) constant-mean-curvature surface in $\mathbf{R}^3$, the famous Wente Torus counterexample to Hopf's Conjecture. Shortly afterwards, N. Kapouleas [K1, K2] provided more examples by a beautiful geometric construction coupled with deep analysis. His construction essentially connects spheres with classical unduloids and nodoids of Delaunay. Our Theorem 3.3 extends his approach from the isotropic Euclidean norm to a "hexagonal" norm for which the isoperimetric Wulff shape is a regular hexagonal prism rather than a sphere, as shown in Figure 1.1. Because the surfaces are polyhedral, the underlying analysis is vastly simplified. On the other hand, because the norm is not smooth, the constant-mean-curvature condition is singular, not even local, and some new methods are required.

An earlier paper ([M1]; see [M5]) provided "cylindrical" unduloids and nodoids. Since the norm and surfaces were rotationally symmetric, the analysis reduced to curves in the plane. Here we must prove equilibrium under a larger class of variations. First we show that the variations can be localized to one face at a time. Second we use an isoperimetric lemma (2.2) to reduce to constant normal translation of the face. Under such variations, the surfaces are in equilibrium by construction. (See the proof of Proposition 3.2.)

**1.1. Acknowledgements.** I would like to thank the very helpful referee for corrections. This work is partially supported by a National Science Foundation grant.

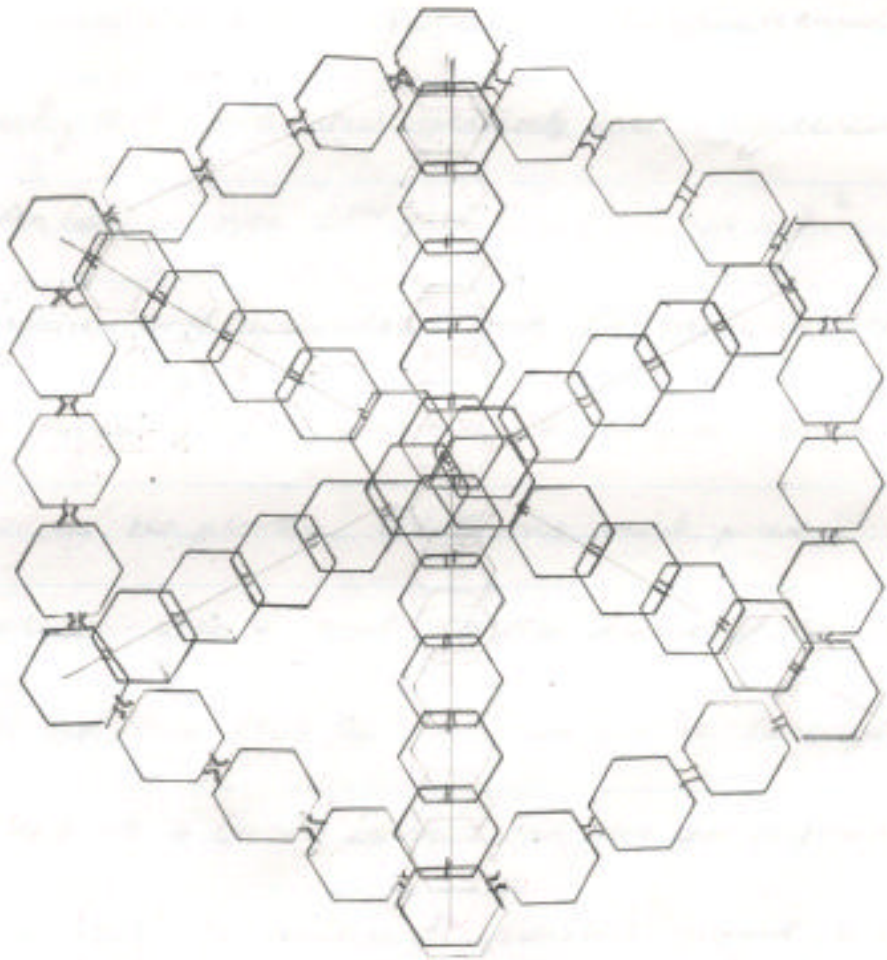

**Figure 1.1**
**A compact immersed surface of hexagonal constant mean curvature. Each hexagon represents a hexagonal prism viewed from above. The outer edges are hexagonal unduloids. The diagonals are hexagonal nodoids, which cross unawares at the center.**

## 2. Hexagonal Norm and First Variation

Section 2 discusses norms and variations, and provides a new useful isoperimetric lemma (2.2). For general background see [M3, Chapt. 10] and [T].

**2.1. Norms and first variation.** Let $\Psi$ be a (continuous) norm on $\mathbf{R}^{n+1}$. An n-dimensional rectifiable surface S with unit normal **n** has associated energy

$$(1) \qquad \Psi(S) = \int_S \Psi(\mathbf{n}) \, dA.$$

The Wulff shape W uniquely provides the least-area way to enclose given volume V. It is the unit ball for the norm dual to $\Psi$.

Following [BdC] and [M3], we say that S has $\Psi$ *mean curvature 1* if the first variation $\Psi^1(\mathbf{v})$ of $\Psi - nV$ is nonnegative for any smooth variation vectorfield **v** of compact support. We say "nonnegative" rather than "zero" because if $\Psi$ is not smooth, $\Psi^1$ is not linear, although it is well defined and homogeneous. For example a small smooth bump **v** on the top of a square Wulff shape in $\mathbf{R}^2$ (with exceptional vertical normal where $\Psi$ is not smooth) will have $\Psi^1(\mathbf{v})$ and $\Psi^1(-\mathbf{v})$ both positive.

Suppose that S is a *planar* piece of surface with normal **n** and piecewise smooth boundary $\partial S$, **u** is the tangential component of **v**, and $v\mathbf{n}$ is the normal component. Differentiation of (1) yields

$$\Psi^1(\mathbf{u} + v\mathbf{n}) = \int_S D\Psi \cdot (-\nabla_S v + (\text{div}_S \mathbf{u})\mathbf{n}) - \int_S n\mathbf{v}.$$

Since $\Psi$ is homogenous, $D\Psi$ is linear on convex combinations of **n** and a tangential vector and

$$(2) \quad \Psi^1(\mathbf{u} + v\mathbf{n}) = \int_S D\Psi \cdot (-\nabla_S v) + \int_S D\Psi \cdot ((\text{div}_S \mathbf{u})\mathbf{n}) - \int_S n\mathbf{v}$$

$$= \int_S D\Psi \cdot (-\nabla_S v) + \int_S \Psi(\mathbf{n})(\text{div}_S \mathbf{u}) - \int_S n\mathbf{v}$$

$$= \int_S D\Psi \cdot (-\nabla_S v) + \int_{\partial S} \Psi(\mathbf{n}) \, \mathbf{u} \cdot \nu - \int_S n\mathbf{v},$$

where $\nu$ is the outward unit conormal. In particular if **w** is tangential,

$$(3) \qquad \Psi^1(\mathbf{v}) = \Psi^1(\mathbf{w}) + \Psi^1(\mathbf{v} - \mathbf{w}).$$

(The situation is more complicated for a *nonplanar* piece of surface.)

The following isoperimetric lemma will be needed to treat general variations in the proof of Proposition 3.2.

**2.2. Lemma.** *Let $\Psi_1$ denote the $l^1$ norm: $\Psi_1(x,y) = |x| + |y|$. For $0 < x_0 < x_1$, $0 < y_0 < y_1$, consider the planar rectangular annulus (Figure 2.1)*

$$\Omega = \{|x| \le x_1, |y| \le y_1\} - \{|x| < x_0, |y| < y_0\}.$$

*There exists an $\varepsilon > 0$ such that if $x_1$ and $y_1$ are between 1 and 100 and $x_0$ and $y_0$ are less than $\varepsilon$, then for all $\Omega' \subset \Omega$,*

(1) $$\frac{\Psi_1(\partial\Omega')}{\text{area } \Omega'} \ge \frac{\Psi_1(\partial\Omega)}{\text{area } \Omega}.$$

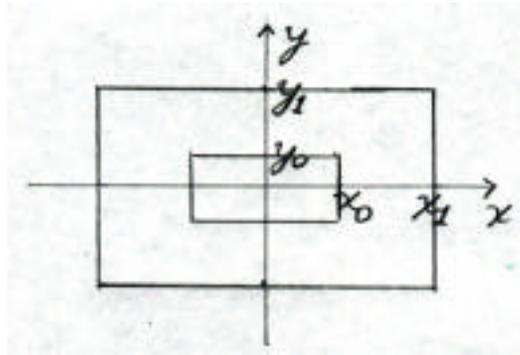

**Figure 2.1**
The planar rectangular annulus

*Proof.* Since two segments meeting at a right angle are $\Psi_1$ minimizing (without area constraint), we may assume that $\partial\Omega'$ consists of vertical and horizontal line segments with no interior concave corners. We may assume that $\Omega'$ is either a square of side at most min$\{x_1-x_0, y_1-y_0\}$, a rectangle of width $x_1-x_0$, a rectangle of height $y_1-y_0$, unions of such, or an annulus $\{|x| \le x_2, |y| \le y_2\} - \{|x| < x_0, |y| < y_0\}$. Since for all such families $\Psi_1$ is a concave function of area, it suffices to check the nontrivial endpoints of each family:

(i)     a square of side say $x_1-x_0$,

(ii)    a rectangle say $x_1-x_0$ by $2y_1$,

(iii)   an L shape with arms of lengths $2x_1$, $2y_1$ and widths $y_1-y_0$, $x_1-x_0$,

(iv)    a U shape with say base $2x_1$, height $2y_1$, and hole $2x_0$ by $2y_1-y_0$.

Since (i) is contained in the family of rectangles, it suffices to check cases ii–iv, where it is sufficient and trivial to check that (1) holds strictly when $x_0 = y_0 = 0$, as in Figure 2.2.

*Remark.* Lemma 2.2 is false in general, e.g. if $x_1 \gg x_0 = y_0 = y_1 = 1$, case (ii).

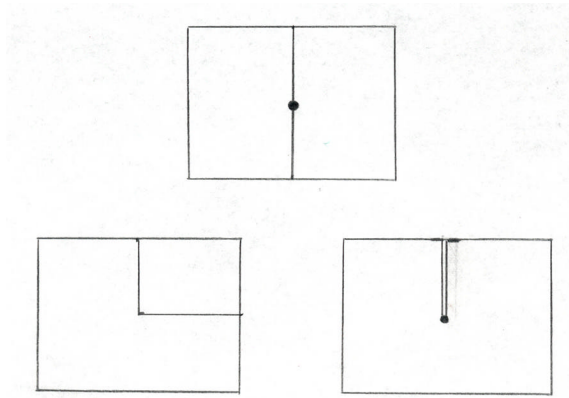

**Figure 2.2**
**Extreme cases of subregions of the rectangle (degenerate rectangular annulus).**

## 3. Hexagonal Surfaces of Delaunay and Kapouleas

Proposition 3.2 provides hexagonal polyhedral analogs of the classical unduloids and nodoids of Delaunay. Following Kapouleas [K1, K2], Theorem 3.3 uses them to construct an immersed compact hexagonal polyhedral surface of constant hexagonal mean curvature.

**3.1. The hexagonal norm.** The *hexagonal norm* on $\mathbf{R}^3$ is characterized by Wulff shape W a right prism of height 2 with a regular hexagonal base of unit inradius and sides $2/\sqrt{3}$ (at $\theta = k\pi/3$). For the unit normal **n** to a face of W, $\Phi(\mathbf{n}) = 1$.

**3.2. Proposition** (Hexagonal surfaces of Delaunay). *For the hexagonal prism Wulff shape, there are one-parameter families of periodic polyhedral surfaces of generalized constant mean curvature 1, the embedded hexagonal unduloids and the immersed hexagonal nodoids of Figure 3.1. As the parameter r goes to 0, they approach a chain of unit hexagonal prisms of height 2 (the Wulff shape).*

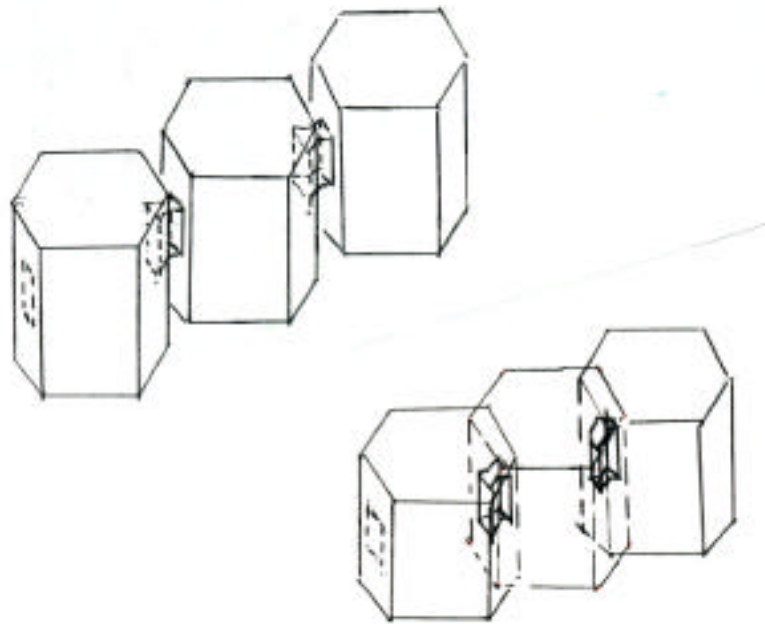

**Figure 3.1**
**Hexagonal unduloid and nodoid of hexagonal constant mean curvature 1.**

*Construction of the unduloids for small r.* **Wulff shapes, with slightly perturbed dimensions, will be connected by polyhedral tubes with parameter r, as shown in Figure 3.2. The affected rectangular face will have width R (slightly less than $2/\sqrt{3}$), height Q (slightly less than 2), and an r×q rectangular hole at the center. The connecting tube will have a flat top and bottom and sides of width s at 120 degrees to the R×Q rectangular faces and to each other. The other large rectangular faces will have width S (slightly greater than $2/\sqrt{3}$).**

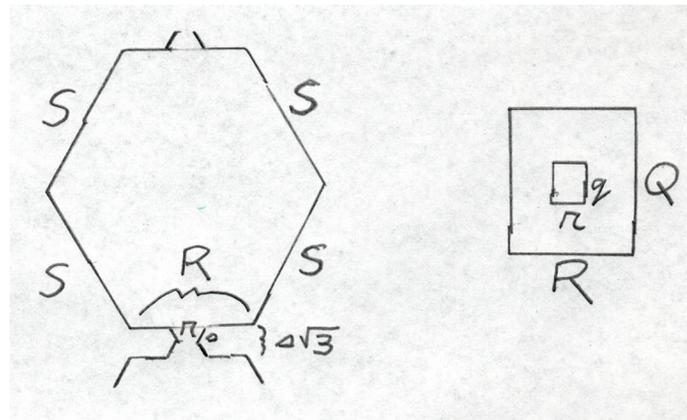

**Figure 3.2**
**To construct hexagonal unduloids, hexagonal prisms (here viewed from above) are connected by small polyhedral tubes of width r, length $s\sqrt{3}$), and height q.**

      **The values of Q, R, S, q, and s as functions of r are chosen to make the first variation of     − 2V vanish for translations of the various faces. Since there are five types of faces, there are five nonlinear equations in the five unknowns. A trivial solution is given by Q = 2, R = S = $2/\sqrt{3}$, q = r = s = 0, and a routine implicit function theorem calculation shows that there are smooth solutions Q < 2, R < $2/\sqrt{3}$, S > $2/\sqrt{3}$, q > 0, s > 0, for small r > 0.**

*Construction of the nodoids for small r.* The nodoids' construction parallels the unduloids', but the Wulff shapes overlap, as in Figures 1.1 and 3.3, resulting in opposite signs in several terms of the formulas, and smooth solutions $Q > 2$, $R > 2/\sqrt{3}$, $S < 2/\sqrt{3}$, $q > 0$, $s > 0$, for small $r > 0$.

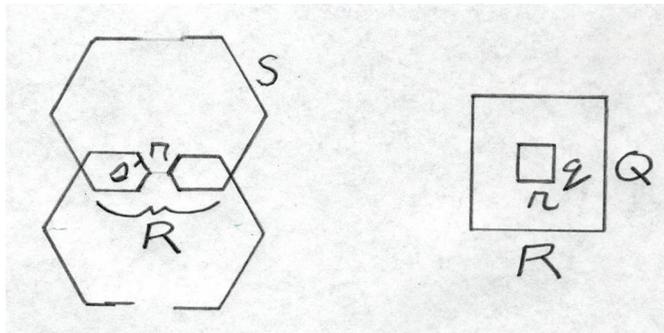

Figure 3.3
To construct hexagonal nodoids, overlapping hexagonal prisms (here viewed from above) are connected by small polyhedral tubes of width r, length $s\sqrt{3}$), and height q.

*Proof of hexagonal mean curvature 1.* We consider a nodoid; an unduloid is similar and easier. On each face $F_i$ choose a smooth nontangential vectorfield $e_i$ which along each edge is parallel to the adjacent face. This means that at a vertex, $e_i$ is parallel to the edge between the two adjacent faces.

Given any variation vectorfield **v**, on each face $F_i$ express **v** as a sum $\mathbf{v} = \mathbf{v}_i + \mathbf{w}_i$ of components in the $e_i$ and tangential directions. By linearity 2.1(3,2), the first variation due to the tangential component can be separated from $\frac{1}{F_i}(\mathbf{v}_i)$ and converted to a boundary integral

(1)
$$\int_{F_i} \mathbf{w}_i \cdot \nu_i$$

( $(\mathbf{n}_i) = 1$).

On the boundary edge $F_i \cap F_j$ ($j \neq i$), let $e_0$ be a vector tangent to that edge, yielding a frame ($e_i$, $e_j$, $e_0$). Express **v** as a sum $\mathbf{v} = \mathbf{v}_i + \mathbf{v}_j + \mathbf{v}_0$ of components in those three directions. Since $e_j$ and $e_0$ are both tangent to $F_i$, the meaning of $\mathbf{v}_i$ has not changed, and $\mathbf{w}_i = \mathbf{v}_j + \mathbf{v}_0$ (and similarly for $\mathbf{v}_j$ and $\mathbf{w}_j$). Since $\mathbf{v}_0$ is perpendicular to $\nu_i$, the contribution from (1) to the $F_i \cap F_j$ portion of the boundary of $F_i$ is

$$\int_{F_i \cap F_j} \mathbf{v}_j \cdot \nu_i .$$

**Thus all contributions to the first variation are expressed in terms of the various $v_i$. Therefore it suffices to prove that for each i, the contribution $\delta^1(v_i)$ to the first variation in terms of $v_i$ is nonnegative:**

$$\delta^1(v_i) = \delta^1_{F_i}(v_i) + \int_{F_i \cap F_j} v_i \cdot v_j \geq 0.$$

**By 2.1(2),**

$$\delta^1(v_i) = \int_{F_i} D\Phi(-\nabla v) + \int_{F_i} v_i \cdot \nu_i + \int_{F_i \cap F_j} v_i \cdot \nu_j - 2\int_{F_i} \Phi v,$$

**where $v = v_i \cdot n_i$ denotes the normal component of $v_i$. By construction, if v is constant, $\delta^1(v_i) = 0$, so by adding a constant to v we may assume that $v > 0$. By scaling, we may assume that $0 < v \leq 1$. We treat the case where $F_i$ is the rectangular annulus $\Omega$; the other cases are similar and easier. We assume the $\Omega$ sits nicely in $\mathbb{R}^2 \times \{0\} \subset \mathbb{R}^3$, with vertical boundary $B_1$ where $v_i \cdot \nu_i = -(1/\sqrt{3})v$ and $v_i \cdot \nu_j = (2/\sqrt{3})v$, and horizontal boundary $B_2$ where $v_i \cdot \nu_i = 0$ and $v_i \cdot \nu_j = v$. Then**

$$\delta^1(v_i) = \int_\Omega D\Phi(-\nabla v) + \int_{B_1}\left(\frac{2}{\sqrt{3}} - \frac{1}{\sqrt{3}}\right)v + \int_{B_2} v - 2\int_\Omega \Phi v.$$

**Since $D\Phi(x,y,0) = \frac{1}{\sqrt{3}}|x| + |y|$, after scaling $\Omega$ by $1/\sqrt{3}$ in the vertical direction,**

$$\delta^1(v_i) = \int_\Omega \Phi_1(-\nabla v) + \int_B v - 2\int_\Omega \Phi v,$$

**where $\Phi_1(x,y) = |x| + |y|$.**

Let $U_t = \{v \geq t\}$ and for now let $\partial U_t$ denote its boundary in the interior of $\Omega$. We apply a version of the (easy, smooth, codimension 1) co-area formula

$$\int_\Omega |\nabla v| = \int_0^1 |\partial U_t|\, dt$$

**or more generally**

(2) $$\int_\Omega |\nabla v|g = \int_0^1 \int_{\partial U_t} g\, dt$$

**[M2, 3.13] for general norms:**

$$\int_\Omega \Phi_1(-\nabla v) = \int_0^1 \Phi_1(\partial U_t)\, dt,$$

which follows from (2) by taking $g = \frac{\delta_1(\nabla v)}{|\nabla v|}$ (or 1 where $\nabla v = 0$), to obtain

$$\delta^1(v_i) = \int_0^1 (\delta_1(\partial U_t) - 2\,\text{area}\,U_t)\,dt,$$

where now $\partial U_t$ denotes the boundary of $U_t$ in the closed domain W, which automatically incorporates the boundary integral over B. By Lemma 2.2, to prove that $\delta^1(v_i) \geq 0$, it suffices to show that

$$\delta_1(\partial\Omega) - 2\,\text{area}\,\Omega \geq 0.$$

But this is the case $v = 1$, when the first variation vanishes by construction.

**3.3. Theorem** (Hexagonal surfaces of Kapouleas). *For the hexagonal prism Wulff shape, there is a compact immersed polyhedral surface of generalized constant mean curvature 1, as in Figure 1.0. The sides are "hexagonal unduloids" while the diagonals are "hexagonal nodoids."*

*Proof.* Six copies of a Wulff shape W, with slightly perturbed dimensions, at the vertices of a big regular hexagon X, will be connected by unduloids on the edges of X and nodoids on the diagonals of X, with parameters $r_1$, $r_2$. The faces meeting the unduloids will have width $R_1$ (slightly less than $2/\sqrt{3}$), the face meeting the nodoid will have width $R_2$ (slightly greater than $2/\sqrt{3}$), and the other three faces will have width $S = (R_1 + R_2)/2$ (slightly greater than $2/\sqrt{3}$). The height $Q_0$ will be slightly less than 2.

The values of $Q_0$, $R_1$, $R_2$, and $r_2$ as functions of $r_1$ are chosen to make the first variation of $\Phi - 2V$ vanish for normal translations of the various faces. There are four types of faces (three types of sides and the top), yielding four equations in four unknowns. A trivial solution is given by $Q_0 = 2$, $R_1 = R_2 = 2/\sqrt{3}$, $r_2 = r_1 = 0$, and a routine implicit function theorem calculation shows that there are smooth solutions $r_2 \sim r_1$, $Q_0 < 2$, $R_1 < 2/\sqrt{3}$, $R_2 > 2/\sqrt{3}$ for small $r_1 > 0$. As in the proof of Proposition 3.2, it follows that our surface has constant mean curvature 1 (for all smooth variations).

The only problem is that as $r_1$ increases from 0, the unduloids get longer, the nodoids get shorter, and they have slightly different distances to traverse. To compensate, we can use for example unduloids and nodoids of many periods, with one extra period in each nodoid, which will fit for some small $r_1$. (Alternatively we could connect them to the anchor prisms W at the vertices of X a bit off center.)